\documentclass[11pt]{article}

\usepackage{latexsym,amssymb}
\usepackage{amsmath}
\usepackage{amscd}
\usepackage{hyperref}

\input xypic
\newcommand{\rar}{\rightarrow}
\newcommand{\lar}{\longrightarrow}

\newcommand{\llar}{-\kern-5pt-\kern-5pt\longrightarrow}
\newcommand{\surjects}{\twoheadrightarrow}

\newtheorem{Theorem}{Theorem}[section]
\newtheorem{Lemma}[Theorem]{Lemma}
\newtheorem{Corollary}[Theorem]{Corollary}
\newtheorem{Proposition}[Theorem]{Proposition}

\newtheorem{Remark}[Theorem]{Remark}
\newtheorem{Example}[Theorem]{Example}
\newtheorem{Definition}[Theorem]{Definition}

\newtheorem{Notation}[Theorem]{Notation}

\def\demo{\noindent{\bf Proof. }}

\def\edeg{\mbox{\rm edeg}}

\def\gr{\mbox{\rm gr}}

\def\height{\mbox{\rm height}}
\def\Hom{\mbox{\rm Hom}}
\def\ker{\mbox{\rm ker}}

\def\m{\mathfrak{m}}

\def\QED{\hfill$\Box$}
\def\qed{\QED}

\def\Rees{\mbox{$\mathcal{R}$}}

\def\Symi{\mbox{\rm Sym}}

\def\AA{{\bf A}}

\def\BB{{\bf B}}

\def\FF{{\bf F}}
\def\HH{{\bf H}}

\def\TT{{\bf T}}
\def\ff{{\bf f}}
\def\pp{{\bf p}}
\def\qq{{\bf q}}
\def\g2{{\bf g}}
\def\hh{{\bf h}}

\def\aa{{\bf a}}
\def\vv{{\bf v}}

\def\zz{{\bf z}}

\begin{document}

\title{{\sc The Equations of Almost Complete Intersections}
\vspace{-1mm}}\footnotetext{2000 AMS {\it Mathematics Subject
Classification}: 13B10, 13D02, 13H10,  13H15, 14E05. }
\vspace{-1mm}

\author{
{\normalsize\sc Jooyoun Hong}
\vspace{-0.75mm}\\
{\small Department of Mathematics}\vspace{-1.4mm} \\
{\small Southern Connecticut State University}\vspace{-1.4mm}\\
{\small 501 Crescent Street}\vspace{-1.4mm}\\
{\small New Haven, CT 06515-1533, U.S.A.}\vspace{-1.4mm}\\
{\small e-mail: {\tt hongj2@southernct.edu}}\vspace{4mm} \and
{\normalsize\sc Aron Simis \thanks{Partially
supported by CNPq, Brazil}}
\vspace{-0.75mm}\\
{\small Departamento de Matem\'atica}\vspace{-1.4mm} \\
{\small Universidade Federal de Pernambuco}\vspace{-1.4mm}\\
{\small 50740-540 Recife, PE, Brazil}\vspace{-1.4mm}\\
{\small e-mail: {\tt aron@dmat.ufpe.br}}\vspace{4mm}
\and
{\normalsize\sc Wolmer V. Vasconcelos\thanks{Partially
 supported by the NSF.}}
\vspace{-0.75mm}\\
{\small Department of Mathematics}\vspace{-1.4mm} \\
{\small Rutgers University}\vspace{-1.4mm}\\
{\small 110 Frelinghuysen Road}\vspace{-1.4mm}\\
{\small Piscataway, NJ 08854-8019, U.S.A.}\vspace{-1.4mm} \\
{\small e-mail: {\tt vasconce@math.rutgers.edu}}\vspace{4mm}
}

\date{June 8, 2009}

\maketitle

\begin{abstract} In this paper we examine the role of four Hilbert
functions in the determination of the defining relations of the Rees
algebra of  almost complete intersections of finite colength. Because
 three of the corresponding modules are Artinian,
some of these  relationships are very effective,
opening up tracks to the determination of
the equations and also to processes of going from homologically
defined sets of equations to higher degrees ones assembled
by resultants.
\end{abstract}
\maketitle

\section{Introduction}

\noindent
Let $R$ be a Noetherian ring  and let $I$ be an
 ideal. By the {\em equations} of $I$ it is meant a free presentation of
 the Rees algebra $R[It]$ of $I$,
\begin{eqnarray}
 0 \rar L \lar S = R[\TT_1, \ldots, \TT_m] \stackrel{\psi}{\lar}
R[It] \rar 0,  \quad \TT_i \mapsto f_it .
\end{eqnarray}
More precisely,  $L$ is the defining ideal of the Rees
algebra of $I$ but we refer to it simply as the {\em ideal of
equations} of $I$.
The ideal $L$ depends on the chosen set of generators of $I$, but all
of its significant cohomological properties, such as the integers that bound the
degrees of minimal generating sets of $L$, are independent of the presentation
$\psi$.
The examination of $L$ is one pathway to the
unveiling of
the properties of $R[It]$.
It codes the syzygies of all powers of $I$, and therefore
 is a carrier of not just algebraic properties of $I$, but of  analytic
ones
as well. It is also a vehicle to understanding geometric properties of
several constructions built out of $I$, particularly of rational maps.

The search for these equations and their use has attracted
considerable interest by a diverse group of researchers. We just
mention some that have  directly influenced this work. One main
source lies in the work of L. Bus\'e, M. Chardin, D. Cox and
 J. P.  Jouanolou  who
have charted, by themselves or with co-workers, numerous roles of
resultants and other elimination techniques to obtain these
equations (\cite{BuCh}, \cite{BuCoDa}, \cite{BuJou}, \cite{Dandrea},
and references therein.) Another important development was given by
A. Kustin, C. Polini and B. Ulrich, who provided a comprehensive
analysis of the equations of ideals (in the binary case), not
just necessarily of almost complete intersections, but still of ideals
whose syzygies   are
almost all linear (\cite{KPU}). Last,
has been the important work of D. Cox not only for its theoretical value to the
understanding of these equations, but for the  role
it has played in bridging the fields of commutative algebra and
of geometric modelling (\cite{CHW}, \cite{Cox08}).

This is a sequel, although not entirely  a continuation of \cite{syl}. It
deals, using some novel methods,
with  questions in higher dimensions that
 were triggered  in that project, but is mainly concerned with the
more general issues of the structure of the Rees algebras of almost
complete intersections.
Our underlying metaphor  here is to focus on distinguished sets of
equations by examining 4 Hilbert functions associated to the ideal
$I$ and to the coefficients of ts syzygies. It brings considerable
effectivity to the methods by developing explicit formulas for some
of the equations.

 There are  natural and technical  reasons to focus on almost intersections. A good
deal of elimination theory is intertwined with birationality
questions. Now, a regular sequence of forms of fixed degree $\geq 2$
never defines a birational map. Thus, the first relevant case is the
next one, namely, that of an almost complete intersection.  Say,
$I\subset R=k[x_1,\ldots,x_d]$ is minimally generated by forms
$a_1,\ldots, a_d, a_{d+1}$ of fixed degree, where $a_1,\ldots, a_d$
form a regular sequence.  If these generators define a birational map
of ${\mathbb P}^{d-1}$ onto its image in ${\mathbb P}^d$ then
any set of forms of this degree containing
$a_1,\ldots, a_d, a_{d+1}$ still defines a birational map onto its
image.  Thus, almost complete intersections give us in some sense the
hard case.

 Note that these almost complete intersections have maximal
codimension, i.e., the ideal $I$ as above is $\m$-primary, where
$\m=(x_1,\ldots,x_d)$. The corresponding rational map with such base
ideal is a regular map with image a hypersurface of ${\mathbb
P}^{d}$. However, one can stretch the theory to one more case,
namely, that of an almost complete intersection of
$I=(a_1,\ldots,a_{d-1}, a_d)$ of submaximal height $d-1$.  Here the
corresponding rational map $\Psi_I:{\mathbb P}^{d-1}\dasharrow {\mathbb
P}^{d-1}$ is only defined off the support $V(I)$ (whose geometric
dimension is $0$) and one can ask when this map is birational---thus
corresponding to the notion of a Cremona transformation of ${\mathbb
P}^{d-1}$. In another paper, we will treat these ideals.

\begin{Notation}\rm To describe the problems treated and the solutions given, we give a
modicum of notation and terminology. When not obvious,
we will point out which characteristics of $k$ to avoid.
\begin{itemize}
\item $R=k[x_1, \ldots, x_d]_{\m}$ with $d \geq 2$ and $\m=(x_1,
\ldots, x_d)$.  \item $I=(a_1, \ldots, a_d, a_{d+1})$ with
$\deg(a_i)=n$ for all $i$, $\height(I)=d$ and $I$ minimally generated
by these forms.  \item Assume that $J=(a_1, \ldots, a_d)$ is a
minimal reduction of $I=(J, a)$, that is $I^{r+1}=JI^r$ for some
natural number $r$.
  \item $ \bigoplus_{j} R(-n_j)
\stackrel{\varphi}{\longrightarrow} R^{d+1}(-n) \longrightarrow I
\longrightarrow 0$ is a free minimal presentation of $I$.  \item $
S=R[\TT_1, \ldots, \TT_{d+1}]$ and $L=\ker(S\surjects R[It])$ via
$T_j\mapsto a_jt$, where $R[It]$ is the Rees algebra of $I$; note that
$L$ is a homogeneous ideal in the standard grading of $S$ with
$S_0=R$.  \item $L_i$ : $R$--module generated by forms of
$L$ of
degree $i$ in $\TT_j$'s. For example, the degree $1$ component of $L$
is the ideal of entries of the matrix product
\begin{eqnarray*}
(L_1)=I_1([\TT_1 \quad \cdots \quad
\TT_{d+1}]\cdot \varphi).\end{eqnarray*}
\item $\nu(L_i)$ denotes the minimal number of {\em fresh} generators
of $L_i$. Thus $\nu(L_2)$ is the minimal number of generators of the
$R$-module $L_2/S_1L_1$.
\item
The {\em elimination degree} of $I$ is
\begin{eqnarray} \label{edeg}
\edeg(I)&=&\inf
\{ i \;|\; L_i \nsubseteq \m S  \}.
\end{eqnarray}
One knows quite generally that $L= (L_1):
\mathfrak{m}^{\infty}$.
 A {\em secondary elimination
degree} is an integer $r$ such that $L= (L_1):
\mathfrak{m}^r$, i.e., an integer at least as large as the stabilizing exponent
of the saturation.
\item  The {\em special fiber} of $I$
is the ring $\mathcal{F}(I)= R[It]\otimes R/\mathfrak{m}$. This is
a hypersurface ring
\begin{eqnarray} \label{elimequ}
 \mathcal{F}(I) = k[T_1, \ldots, T_{d+1}]/(\ff(\TT)),\end{eqnarray}
 where
$\ff(\TT)$ is an irreducible polynomial of degree $\edeg(I)$.
Then $\ff(\TT)$ is called the {\em elimination equation} of $I$, and may
be taken as  an element of $L$.
\end{itemize}
\end{Notation}

Besides the syzygies of $I$, $\ff(\TT)$ may be considered the most
significant of the {\em equations} of $I$. Determining it, or at
least its degree, is one of the main goals of this paper.
The enablers, in our treatment, are   four Hilbert functions associated
to $I$, those of $R/I$, $R/J:a$, $R/I_1(\varphi)$ and the
Hilbert-Samuel function defined by $I$. Each encodes, singly or in
conjunction, a  different
aspect of $L$.

\medskip

 In order to
describe the other relationships between the invariants of the ideal $I$ and
its equations, we make use the following diagram:

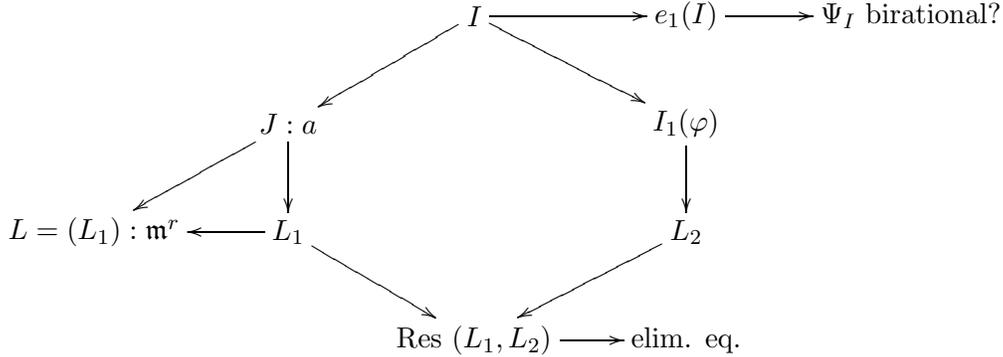
\[\diagram
& & I \dlto \rto \drto & e_1(I)\rto & \mbox{\rm $\Psi_I$ birational?} \\
& J:a \dlto \dto & & I_1(\varphi)\dto & \\
L=(L_1): \mathfrak{m}^r &  L_1\lto \drto & & L_2 \dlto  & \\
& & \mbox{\rm Res }(L_1,L_2)\rto & \mbox{\rm elim. eq.} & \\
\enddiagram
\]

At the outset and throughout there is the role
 played by
 the Hilbert function of $J:a$, which besides that of directing all the
 syzygies of $I$,  is the encoding  of an integer $r$
such that $L= (L_1):\mathfrak{m}^r$. According to
Theorem~\ref{secondedeg}, $r$ can be taken as  $r=\epsilon+1$,
 where $\epsilon$  is the {\em socle degree} of
$R/J:a $. Since $L$ is expressed as the quotient of two
Cohen-Macaulay ideals, this formula  has shown in practice to be an effective
tool to determine  saturation.

A persistent  question is that of how to obtain higher degree
generators from the syzygies of $I$.
We will provide an iterative approach to generate the successive
components of $L$:
\[ L_1 \leadsto L_2 \leadsto L_3 \leadsto \cdots .\]
This is not an effective process, except for the step
$ L_1\leadsto L_2$.
 The more interesting development  seems to be the use
of the syzygetic lemmas to
obtain $\delta(I)=L_2/S_1L_1$ in the case of almost complete
intersections. This is a formulation of the {\em method of moving
quadrics} of several authors. The novelty here is the use of the
Hilbert function of the
ideal $I_1(\varphi)$ to understand and conveniently express
$\delta(I)$. Such  level of detail was not present even in earlier
treatments of $\delta(I)$.
The {\em syzygetic lemmas}   are observations based on the Hilbert
functions of $I$ and of $I_1(\varphi)$ to allow a description of
$L_2$. It converts the expression (\ref{deltaI})
\begin{eqnarray} \delta(I) = \Hom_{R/I}(R/I_1(\varphi),
H_1(I)),\end{eqnarray}
where $H_1(I)$ is the canonical module of $R/I$ (given by the
syzygies of $I$), into a set of generators of $L_2/S_1L_1$.
It is fairly effective in the case of binary forms, many cases of  ternary and some
quaternary forms, as we
shall see. In these cases, out of $L_1$ and $L_2$ we will be able to
write the elimination equation in the form of a resultant
\begin{eqnarray}
\textrm{\rm Res }(L_1,L_2),\end{eqnarray} or as one of its factors.
 We  prove the
non-vanishing of this determinant under three different situations.
In the case of binary ideals, whose syzygies are of arbitrary
degrees, we give a far-reaching generalization of \cite{syl}. Here we
make use of one of the distinguished submodules of $\delta(I)$,
\[ \delta_s(I) = \Hom_R(R/\mathfrak{m}^s, \delta(I))\hookrightarrow
\delta(I),\]
where  $s$ is the order of the ideal $I_1(\varphi)$. While
$\delta(I)$ accounts for the whole of $L_2$, $\delta_s(I)$ collects
the forms to be assembled in a resultant.
 For
an ideal $I$ of $k[x,y]$, generated by $3$ forms of degree $n$,
we build out of
$L_1$ and $L_2$, in a straightforward manner,  a nonzero
polynomial of $k[\TT_1, \TT_2, \TT_3]$ in $L$, of degree $n$ (Theorem~\ref{detB}).

The other results, in higher dimension,
  require that the content of syzygies ideal $I_1(\varphi)$ be
a power of the maximal ideal $\mathfrak{m}$, that is
$\delta_s(I)=\delta(I)$.  Thus, in the case of
ternary forms of degree $n$ whose content ideal $I_1(\varphi)=
\mathfrak{m}^{n-1}$, our main result (Theorem~\ref{detB3s}) proves the
non-vanishing of $\mbox{\rm Res }(L_1,L_2)$ without appealing to
conditions of genericity (but with degree constraints). A final result
is  very special to quaternary forms (Theorem~\ref{detB42}).

\section{Syzygies and Hilbert Functions}

Before engaging in the above questions proper, we will outline the
basic homological and arithmetical data involved in these
ideals.

\subsubsection*{The resolution}

The general format of the resolution of $I$ goes as follows. First,
note that $J:a$ is a Gorenstein ideal, and that $I=J: (J:a)$. Since
Gorenstein ideals, at least in low dimensions, have an amenable
structure, it may be desirable to look at $J:a$
as a building block to $I$  and its equations. In this arrangement
the syzygies of $I$ will be organized in terms of those of $J$ and of
$J:a$.

\medskip

Let us recall how the resolution of $I$ arises as a mapping cone of
the Koszul complex $\mathbb{K}(J)$ and a minimal resolution of
$J:a$ (this was first given in \cite{pszpiro}; see also
\cite[Theorem A.139]{compu}).

\begin{Theorem} \label{Dubreil}
 Let $R$ be a Gorenstein local ring, let ${\mathfrak a}$ be a
perfect ideal of height $g$ and let $\mathbb{F}$ be a minimal free
resolution of $R/{\mathfrak a}$. Let ${\bf z}= z_1, \ldots, z_g \subset {\mathfrak a}$ be a
regular sequence, let
$\mathbb{K}=\mathbb{K}({\bf z}; R)$ be the corresponding Koszul complex,
 and let $u\colon  \mathbb{K} \rar \mathbb{F}$ be the
comparison mapping induced by the inclusion $({\bf z}) \subset I$.
Then the dual $\mathbb{C}(u^*)[-g]$ of the mapping cone of $u$, modulo the
subcomplex $u_0\colon R \rar R$, is a free resolution of length $g$ of
$R/({\bf z})\colon {\mathfrak a}$.
Moreover, the canonical module of
$R/({\bf z})\colon {\mathfrak a}$ (modulo shift in the graded case)
is ${\mathfrak a}/({\bf z})$.
\end{Theorem}

\begin{Remark} \label{Dubrmk} {\rm For later reference,
we point out
three observations when  the ideal
is the above almost complete intersection $I$.
\begin{enumerate}

\item If $J:a = (b_1, \ldots, b_d)$, writing
\[ [a_1, \ldots, a_d] = [b_1, \ldots, b_d]\cdot \phi, \]
gives
\[ I = (J, \det(\phi)),\]
so that $I$ is a Northcott ideal.

\item If $I$ is a generated by forms of degree $n$, a choice for $J$
is simply a set of forms $a_1, \ldots, a_d$ of degree $n$ generating
a regular sequence. Many of the features of
$I$--such as the ideal $I_1(\varphi)$--can be read off $J:a$.
Namely, $I_1(\varphi)$ is the sum of $J$ and  the coefficients arising
in the expressions of $a(J:a)\subset J$.

\item Since $J:a$ is a  Gorenstein ideal, its rich structure in
dimension $\leq 3$ is fundamental to the study in these cases.

\end{enumerate}

}\end{Remark}

\subsubsection*{Hilbert functions}

There are four Hilbert functions related to the ideal $I$ that are
significant in this paper, the first three of the Artinian modules
$R/I$, $R/I_1(\varphi)$ and $R/J:a$, and are  therefore of easy
manipulation. Their interactions will be a mainstay of the paper.

\medskip

 The first elementary observation, whose proof is well-known as to be omitted,
 will be useful when we need the Hilbert function of the
canonical module of $R/I$.

\begin{Proposition} \label{Dubrmk2}
If $I$
\ is generated by forms of degree $n$ of $k[x_1,
\ldots, x_d]$, the Hilbert function of $R/I$ satisfies
\begin{eqnarray*}  H_{R/I}(t) & = &  H_{R/J}(t) - H_{I/J}(t) \\
&=& H_{R/J}(t) - H_{R/J:a}(t-n).
\end{eqnarray*}
\end{Proposition}

The fourth Hilbert function is that of the associated graded ring
\[ \gr_I(R) = \bigoplus_{m\geq 0} I^m/I^{m+1}.\]
It affords the Hilbert--Samuel polynomial ($m\gg 0$)

\[ \lambda(R/I^{m+1}) = e_0(I){{d+m}\choose{d}} - e_1(I)
{{d+m-1}\choose{d-1}} + \textrm{lower degree terms of $m$},\]
where $e_0(I)$ is the
multiplicity of the ideal $I$. A related Hilbert polynomial is that
associated to the integral closure filtration $\{\overline{I^m}\}$:
 \[ \lambda(R/\overline{I^{m+1}}) =
 \overline{e}_0(I){{d+m}\choose{d}} - \overline{e}_1(I)
{{d+m-1}\choose{d-1}} + \textrm{lower degree terms of $m$}.\]
For an $\mathfrak{m}$-primary ideal $I$ generated by forms of degree $n$, $\overline{I^m}=
\mathfrak{m}^{nm}$, so the latter coefficients are really invariants of the ideal
$\mathfrak{m}^{n}$ and one has
\begin{eqnarray*}
e_0(I) &=& \overline{e}_0(I) = n^d \\
e_1(I) &\leq & \overline{e}_1(I) =   {\frac{d-1}{2}}(n^d-n^{d-1}).\\
\end{eqnarray*}

The case of equality $e_1(I)=\overline{e}_1(I)$ has a straightforward
(and general) interpretation in terms of the corresponding Rees
algebras.

\begin{Proposition} \label{r1} Let $(R, \mathfrak{m})$ be an
analytically unramified normal local domain of dimension $d$, and
 let $I$ be an
$\mathfrak{m}$-primary ideal. Then $e_1(I)=\overline{e}_1(I)$ if and
only if $\Rees(I)$ satisfies the condition $(R_1)$ of Serre.
\end{Proposition}

\demo Let $\AA= \Rees(I)$, set $\BB$ for its integral closure. Then $\BB$
is a finitely generated $\AA$-module. Consider the exact sequence
\[ 0 \rar \AA \lar \BB \lar C = \BB/\AA \rar 0. \]
Since $e_0(I)= \overline{e}_0(I)$, $C$ is a graded $\AA$-module of
dimension $\leq d$. If $\dim C=d$, its multiplicity $e_0(C)$ is derived from
the Hilbert polynomials above as $e_0(C) = \overline{e}_1(I)-e_1(I)$.
This sets up the assertion since $\AA$ and $\BB$ are equal in
codimension one if and only if $\dim C< d$. \QED

\bigskip

This permits stating \cite[Proposition 3.3]{syl} as follows (see also
\cite{ehrhart} for degrees formulas).

\begin{Proposition}
Let $R=k[x_1, \ldots, x_d]$ and let $I = (f_1,
\ldots, f_{d+1})$ be an ideal of finite colength, generated by forms
of degree $n$. Denote by $\mathcal{F}$ and $\mathcal{F}'$ the special
fibers of $\Rees(I)$ and $\Rees(\mathfrak{m}^n)$, respectively. The following
conditions are equivalent:
\begin{enumerate}
\item[{\rm (i)}] $[\mathcal{F}':\mathcal{F}]=1$, that is, the rational mapping
\[ \Psi_I = [f_1: f_2: \cdots : f_{d+1}] : \mathbb{P}^{d-1}\dasharrow
\mathbb{P}^d\]
is birational onto its image$\,${\rm ;}

\item[{\rm (ii)}] $\deg \mathcal{F}= n^{d-1}\,${\rm ;}

\item[{\rm (iii)}] $e_1(I) =  {\frac{d-1}{2}}(n^d-n^{d-1})\,${\rm ;}

\item[{\rm (iv)}] $\Rees(I)$ is non-singular in codimension one.
\end{enumerate}
\end{Proposition}

A great deal of this investigation is to determine $\deg
\mathcal{F}$, which as we referred to earlier is the {\em elimination degree} of
$I$ (in notation, edeg$(I)$). Very often this turns into
explicit formulas for the elimination equation.

\subsubsection*{Secondary elimination degrees}

A solution to some of questions raised above resides in the
understanding of the primary decomposition of $(L_1)$, the defining
ideal of the symmetric algebra ${\rm Sym}(I)$ of $I$ over $S=R[{\bf T}]$.
As pointed out earlier, we know that $L=(L_1):\mathfrak{m}^{\infty}$
as quite generally a power of $I$ annihilates $L/(L_1)$.
The $L$-primary component is therefore $L$ itself and $(L_1)$ has
only two primary components, the other being $\mathfrak{m}S$-primary.
On the other hand, according to
\cite[Proposition 2.2]{syl}, $(L_1)$ is Cohen-Macaulay so its
primary components are of the same dimension.
Write
\[ (L_1)= L\cap Q,\]
where $Q$ is $\mathfrak{m}S$-primary. It allows for the expression of
$L$ as a saturation of $(L_1)$ in many ways. For example,
 for nonzero $\alpha \in I$ or $\alpha \in
I_1(\varphi)$, we have $L=(L_1) :\alpha^{\infty}$.

\medskip

We give now an explicit saturation by exhibiting integers $r$
 such that $L=
(L_1):\mathfrak{m}^r$, which as we referred to earlier are {\em secondary
elimination degrees}.
Finding its least value is one of our goals in
individual cases.
In the actual practice we have found the
computation effective, perhaps because $L$ is given as the quotient
of two Cohen-Macaulay ideals, the second generated by monomials.

\begin{Theorem}\label{secondedeg}
Let $R=k[x_1, \ldots, x_d]$ and let $I = (f_1,
\ldots, f_{d+1})$ be an ideal of finite colength, generated by forms
of degree $n$. Suppose $J=(f_1, \ldots, f_{d})$ is a minimal
reduction, and set $a=f_{d+1}$. Let $\epsilon$ be the socle degree of
$R/(J:a)$, that is the largest integer $m$ such that $(R/(J:a))_m\neq 0$.
If $r=\epsilon+1$, then
\[ L= (L_1): \mathfrak{m}^r.\]
More precisely,  if some
form $\ff$ of $L_i$ has
coefficients in $\mathfrak{m}^r$, then $\ff\in (L_1)$.
\end{Theorem}

\demo The assumption on $r$ means that $(J:a)_m=({\mathfrak{m}^m})_m$ for
$m\geq r$, that is if $J:a$ has initial degree $d'$ then
\[ \mathfrak{m}^{r} = \sum_{i\geq d'}(J:a)_i \mathfrak{m}^{r-i}.\]

Now any element $\pp\in L$ can be written as
\[ \pp = \hh_p \TT_{d+1}^p + \hh_{p-1}\TT_{d+1}^{p-1} + \cdots + \hh_0, \]
where the $\hh_i$ are polynomials in $\TT_1, \ldots, \TT_{d}$.
If $\hh_1,\ldots, \hh_{p}$ all happen to vanish, then $\hh_0\in (L_1)$ since $f_1, \ldots,
f_{d}$ is a regular sequence.
The proof will consist in reducing to this situation.

To wit, let $u$ be a form of degree $r$ in $ \mathfrak{m}^r$. To prove that
$u\pp\in (L_1)$, we may assume that $u=v\alpha$, with $v\in
\mathfrak{m}^{r-i}$ and $\alpha$ a minimal generator in $(J:a)_i$.

  We are going to replace $v\alpha \pp$ by an
equivalent element of $L$, but of lower degree in $\TT_{d+1}$.
Since $\alpha\in (J:a)$, there is an induced form $\g2\in L_1$
\[ \g2 = \alpha \TT_{d+1} + \hh, \quad \mbox{\rm $\hh$ linear form in $\TT_1, \ldots, \TT_{d}$}.\]
Upon substituting $\alpha \TT_{d+1}$ by $\g2-\hh$, we get an equivalent
form
\[ \qq = \g2_{p-1} \TT_{d+1}^{p-1} + \g2_{p-2}\TT_{d+1}^{p-2} + \cdots + \g2_0, \]
where the coefficients of the $\g2_i$ all lie in $\mathfrak{m}^r$.
Further reduction of the individual terms of $\qq$ will eventually
lead to a form only in the $\TT_1, \ldots, \TT_{d}$.

The last assertion just reflects the nature of the proof.
\QED

\begin{Remark}{\rm An a priori bound for the smallest secondary elemination
degree arises as follows. Let $R=k[x_1, \ldots, x_d]$,  and $I=
(J,a)$ as above. Since $J$ is generated by $d$ forms $f_1, \ldots, f_d$ of degree $n$,
the socle of $R/J$ is determined by the Jacobian of the $f_i$ which
has degree $d(n-1)$. Now from the exact sequence
\[0 \rar (J:a)/J \lar R/J \lar R/(J:a)\rar 0,\]
the socle degree of $R/(J:a)$ is smaller than $d(n-1)$, and therefore
this value gives the bound. In fact, all the examples we examined had $\epsilon+1 = \min\{\; i \mid L=(L_1):\m^{i} \}$, where $\epsilon$ is the socle degree of $R/(J:a)$.
}\end{Remark}

\begin{Example}{\rm Let  $R=k[x_1,x_2,x_3, x_4]$ and
 $\mathfrak{m}=(x_1, x_2, x_3, x_4)$.

\noindent Let $I=(x_1^3,\;\; x_2^3,\;\; x_3^3,\;\; x_4^3,\;\;
x_1^2x_2+x_3^2x_4)$, $J=(x_1^3,\;\; x_2^3,\;\; x_3^3,\;\; x_4^3)$, and $a=x_1^2x_2+x_3^2x_4$. Then
\begin{enumerate}
\item[$\bullet$] Hilbert series of $I$ : \, $1+
4t+10t^2+15t^3+15t^4+7t^5+t^6$.
\item[$\bullet$]  Hilbert series of
$(J:a)$ : \, $1+ 4t+9t^2+9t^3+4t^4+t^5$.
\item[$\bullet$]
Hilbert series of $I_1(\varphi)$ : \, $1+ 4t+7t^2$.
\end{enumerate}
A run with {\it Macaulay2} showed:
\begin{enumerate}
\item $L=(L_1):\mathfrak{m}^6$, exactly as predicted from the Hilbert
function of $J:a$;
\item The calculation yielded $\edeg(I)=9$; in particular, $\Psi_I$ is not
birational.

\end{enumerate}
}
\end{Example}

\subsubsection*{The syzygetic lemmas}

The following material complements and refines some  known facts
(see  \cite{HSV1}, \cite{SV1}, \cite[Chapter 2]{alt}).
Its main purpose is an application to almost complete
intersections.
Since its contents deal with arbitrary ideals, we will momentarily change notation.
Let $I\subset R$ be an ideal generated by $n$ elements $a_1,\ldots, a_n$.
Consider a free presentation of $I$
\begin{equation}\label{presentation}
R^m\stackrel{\varphi}{\lar} R^n\lar I\rar 0
\end{equation}
and let $(L_1)\subset L=\bigoplus_{d\geq 0} L_d \subset S=R[\TT]=R[\TT_1, \ldots,\TT_n]$ denote as
before the presentation ideals of the symmetric algebra
and of the Rees algebra of $I$, respectively, corresponding to the chosen presentation.

\medskip

A starting point is the following
  observation. Suppose $\ff(\TT)= \ff(\TT_1, \ldots,
\TT_n)\in L_d$; write it as
\[ \ff(\TT_1, \ldots, \TT_n) = \ff_1(\TT) \TT_1 + \cdots +\ff_n(\TT)
\TT_n,\]
where $\ff_i(\TT)$ is a form of degree $d-1$.

\medskip

Evaluating $\TT$ at the vector $\aa= (a_1, \ldots, a_n)$ gives a
syzygy of $\aa$ \[ z=(\ff_1(\aa), \ldots, \ff_n(\aa)) \in Z_1,\] the
module of syzygies of $I$, \[ z\in Z_1 \cap I^{d-1}R^n,\] that is,
$z$ is a syzygy with coefficients in $I^{d-1}$.  Note that \[
\widehat{\ff}(\TT)=a_1\ff_1(\TT) + \cdots + a_n \ff_n(\TT) \in
L_{d-1}\cap I\cdot S_{d-1}.\] Conversely, any form ${\hh}(\TT)$ in
$L_{d-1}\cap I\cdot S_{d-1}$ can be lifted to a form $\HH(\TT)$ in
$L_d $ with $\widehat{\HH}(\TT)=\hh(\TT)$.

\medskip

 Such  maps are referred to \index{downgrading and upgrading maps}
as {\em downgrading} and {\em
upgrading}, although they are not always well-defined. However, in some case
it opens the opportunity to calculate some of the higher $L_d$.
\medskip

Here is a useful observation.

\begin{Lemma} \label{syzlemma1} Let
 $\ff(\TT)\in L_d$ and write
\[  \ff(\TT) = \ff_1(\TT) \TT_1 + \cdots + \ff_n(\TT) \TT_n.\]
If
 \[ \ff_1(\TT) a_1 + \cdots +\ff_n(\TT) a_n=0,\]
then $\ff(\TT) \in S_{d-1}L_1$.
\end{Lemma}

\demo The assumption is that $\vv=(\ff_1(\TT), \ldots, \ff_n(\TT))$
is a syzygy of $a_1, \ldots, a_n$ over the ring $R[\TT]$. By
flatness, \[ \vv = \sum_{j} \hh_j(\TT) \zz_j,\] where $\hh_j(\TT)$
are forms of degree $d-1$ and $\zz_j\in Z_1$. Setting $\zz_j =
(z_{1j}, \ldots, z_{nj})$, \[ \ff_i(\TT) = \sum_{j} \hh_j(\TT)
z_{ij},\] and therefore
\begin{eqnarray*}
  \ff(\TT) &=& \ff_1(\TT) \TT_1 + \cdots +\ff_n(\TT) \TT_n\\
&=& \sum_i (\sum_j \hh_j(\TT) z_{ij})\TT_i \\
&=& \sum_{j}\hh_j(\TT) (\sum_{i} z_{ij}\TT_i) \in S_{d-1}L_1.\\
\end{eqnarray*}
\QED

\begin{Corollary} Let
$\hh_j(\TT)$, $1\leq j\leq m$, be a  set of generators
of  $L_{d-1}\cap
IS_{d-1}$. For each $j$, choose a form $\FF_j(\TT)\in L_d$ such that
under one of the operations above $\widehat{\FF}(\TT)=\hh_j(\TT)$. Then
\[ L_d = (\FF_1(\TT), \ldots, \FF_{m}(\TT), L_1S_{d-1}).\]
\end{Corollary}

\demo  For $\ff(\TT)\in L_d$, $\ff(\TT)=
\sum_{i}\TT_i\ff_i(\TT)$, write
\[  \sum_{i}a_i\ff_i(\TT) = \sum_j c_{j}\hh_j(\TT).\]
Applying Lemma~\ref{syzlemma1} to the   polynomial
\[ \ff(\TT) - \sum_{j}c_j\FF_j(\TT)\]
will give the desired assertion. \QED

\bigskip

\noindent This leads to the iterative procedure to find the equations  $L=(L_1, L_2,
\ldots)$ of the ideal $I$.

\bigskip

Let $Z_1=\ker(\varphi)\subset R^n$, where $\varphi$ is as in (\ref{presentation})
and let $B_1$ denote the submodule of $Z_1$ whose elements come from the Koszul relations
of the given set of generators of $I$.
The $R$-module
\[ \delta(I)=  Z_1\cap IR^n/B_1 \]
has been introduced in \cite{Si} in order to understand the Koszul homology with
coefficients in $I$. It is independent of the free presentation of $I$
and as such it has been dubbed the {\em syzygetic module} of $I$.
 The following basic result has been proved in \cite{SV1}.

\begin{Lemma}{\rm (\cite[1.2]{SV1}) (The syzygetic lemma)}  \label{syzlemma} Let $I$ be an ideal with
presentation as above.  Then
\begin{equation}\label{main_syzygetic}
\delta(I)\, \stackrel{\phi}{\simeq} \, {L}_2/ {L}_1S_1.\\
\end{equation}
\end{Lemma}
The mapping $\phi$ is given as follows: For $z= \sum \alpha_i\TT_i$,
$(\alpha_1, \ldots, \alpha_n)\in Z_1\cap IR^n$,
$\alpha_i=\sum_{j=1}^n c_{ij}a_j$,
\[ \phi([z]) = \sum_{i,j}c_{ij}\TT_i\TT_j \in
 {L}_2/S_1 {L}_1.\]
 In particular, ${L}_2/S_1L_1$ is also independent of the free presentation of $I$ and
\[ \nu( {L}_2/S_1 {L}_1)= \nu(\delta(I)).\]

\medskip

We refer to the process of writing the $\alpha_i$ as linear
combination of the $a_j$ as the {\em extraction} of $I$. The
knowledge of the degrees of the $c_{ij}$ is controlled by the
degrees of $\delta(I)$. Note that, quite generally, the kernel of the natural surjection
$\Symi(I) \lar \Rees(I)$
in degree $d$ is $L_d/L_1S_{d-1}$.
However, a more iterative form of bookkeeping of $L$ is through the modules $L_d/S_1L_{d-1}$
that represent the {\em fresh} generators in degree
$d$. Unfortunately, except for the case $d=2$, one knows no explicit expressions for
these modules, hence the urge for different methods to approach the problem.

\subsubsection*{Almost complete intersections}
We now go back to the particular setup of almost complete intersections.
As before, $(R, \mathfrak{m})$ denotes the standard graded polynomial ring
$k[x_1,\ldots, x_d]$ and its irrelevant ideal and $I=(a_1, \ldots,
a_d, a_{d+1})$ is an $\mathfrak{m}$-primary ideal minimally generated
by $d+1$ forms of the same degree. We assume that $J=(a_1, \ldots, a_d)$ is a minimal
reduction of $I$; set $a=a_{d+1}$.

\medskip

Considerable  numerical information about $ {L}_2$ is readily available
in this setup.

\begin{Proposition} \label{syzlemma0} Let $I$ be as above  and
let $\varphi$ be a minimal presentation map as in {\rm (\ref{presentation})}. Then
\[ \nu( {L}_2/S_1 {L}_1)= \nu (I_1(\varphi):
\mathfrak{m})/I_1(\varphi)).\] In particular, if
$I_1(\varphi)=\mathfrak{m}^{s}$, $s\geq 1$, then
\[ \nu( {L}_2/S_1 {L}_1)= {{d+s-2}\choose{d-1}}.\]
Moreover, $\delta(I)$ is generated by the last $s$ graded components of the
first Koszul homology module $H_1(I)$.
\end{Proposition}

\demo Consider the so-called syzygetic sequence of $I$
\[ 0 \rar \delta(I)\lar H_1(I) \lar (R/I)^{d+1} \lar I/I^2 \rar 0. \]
Note that  $H_1(I)$ is isomorphic to the canonical module of $R/I$. Dualizing with
respect to $H_1(I)$ gives the exact complex
\[ H_1(I)^{d+1} \lar \Hom_{R/I}(H_1(I),H_1(I))\simeq  R/I \lar
\Hom_{R/I}(\delta(I), H_1(I)) \rar 0.\] The image in $R/I$ is the ideal
generated by $I_1(\varphi)$, and since $I\subset I_1(\varphi)$,  one has
\begin{eqnarray} \label{deltaI}
 \delta(I) \simeq \Hom_{R/I}(R/I_1(\varphi), H_1(I)).
\end{eqnarray}
It follows that $\delta(I)$ is isomorphic to the canonical module of
$R/I_1(\varphi)$,
and therefore $\nu(\delta(I))$ is the
Cohen-Macaulay type
of $R/I_1(\varphi)$.

In case $I_1(\varphi)=\mathfrak{m}^{s}$,
$\Hom_{R/I}(R/\mathfrak{m}^s, H_1(I))$ cannot have a nonzero element
$u$ in  higher degree as otherwise $\mathfrak{m}^{s-1}u$ would lie in
the socle of $H_1(I)$, a contradiction.
 It follows that $\delta(I)$  is generated
by ${d+s-2}\choose{d-1}$ elements.
\QED

\bigskip

To help identify the generators of $\delta(I)$
requires information about the Hilbert function of $H_1(I)$.
For reference  we use
 the socle degree of
$R/I_1(\varphi)$, which we denote by $p$.
 We
recall that if $(1, d, a_2, \ldots, a_r)$ is the Hilbert function of
$R/I$, that of $H_1(I)$ is $(a_r, \ldots, d, 1)$, together with an
appropriate
shift. Since $\delta(I)$ is a graded submodule of $H_1(I)$, it is convenient to
organize a table as follows:
\begin{eqnarray*}
H_1(I) & : & (a_r, \ldots, a_s, \ldots, a_2, d, 1) \\
\delta(I) &: & (b_p, \ldots, b_1,1),
\end{eqnarray*}
where $b_i\leq a_i$. Note that the degrees are increasing. For
example, the degree of the $r$th component of $H_1(I)$ is the initial
degree $t$ of the ideal $J:a$, while the degree of the $p$th component of
$\delta(I)$ is $t+r-p$.

\subsubsection*{Balanced ideals}

Let us introduce the following concept for easy reference:

\begin{Definition}\rm
An $\mathfrak{m}$-primary ideal $I\subset R$ minimally generated
by $d+1$ elements of the same degree is $s$-{\em balanced} if
 $I_1(\varphi)=\m^s$, where $\varphi$ is the matrix of syzygies of
 $I$.
\end{Definition}

Note  that, due to the Koszul relations,  $s$ is at most the common
degrees of the generators of $I$.
The basic result about these ideals goes as follows.

\begin{Theorem} \label{socdegree} Let $R=k[x_1, \ldots, x_d]$ and
 $I$ an almost complete  intersection of  finite colength
 generated by forms of degree $n$. If $I_1(\varphi)\subset
 \mathfrak{m}^s$, for $s$ as large as possible,
  then:
 \begin{enumerate}
 \item[{\rm (i)}]  The socle degree of $H_1(I)$ is $d(n-1)$;
 \item[{\rm (ii)}]  $\mathfrak{m}^{d(n-1)-s+1} = I
 \mathfrak{m}^{(d-1)(n-1)-s}$;
 \item[{\rm (iii)}] Suppose $I$ is $s$-balanced. Let $r(I)$ be the degree the coefficients of $L_2$.
 Then
 \[ r(I)=(d-1)(n-1)-s, \;\; \mbox{\rm and} \;\; R_{n+r(I)}= I_{n+r(I)}.\]

\end{enumerate}
\end{Theorem}
\demo
(i)  Let $J$ be the minimal reduction of $I$ defined earlier.
The socle degree of $H_1(I)$ is determined from the natural
embedding
\begin{eqnarray*}
H_1(I) \simeq J:a/J \hookrightarrow R/J,
\end{eqnarray*} where the socle of $R/J$, which is also the socle of
any of its nonzero submodules,
 is defined by the Jacobian
determinant
of $d$ forms of degree $n$.

\medskip

(ii) We write $H_1(I)$ and $\delta(I)$ as graded modules (set
$\epsilon=d(n-1)$)
\begin{eqnarray*}
H_1(I) & = &  h_s \oplus h_{s+1} \oplus \cdots \oplus h_{\epsilon} \\
\delta(I) &= & f_{\epsilon -s+1} \oplus \cdots \oplus f_{\epsilon},
\end{eqnarray*} dictated by the fact that the two modules share the
same socle, $h_{\epsilon}=f_{\epsilon}$. One has
$\mathfrak{m}^{\epsilon-s+1} H_1(I) = 0$, hence
$ \mathfrak{m}^{\epsilon-s+1} R/I = 0$,
or equivalently, \[\mathfrak{m}^{\epsilon-s+1}=
I\mathfrak{m}^{\epsilon-n-s+1}.\]

(iii) The degree $r(I)$ of the coefficients of $L_2$ is obtained from the
elements of $f_{\epsilon -s+1}$, and writing them as syzygies with
coefficients in $I$, that is
\[ r(I)= \epsilon-s+1-n = (d-1)(n-1)-s.\]
The last assertion follows from (ii).
\qed

\bigskip

Let us give some consequences of this analysis which will be used later.

\begin{Corollary} \label{wellbal} Let $R=k[x_1, \ldots, x_d]$ be a ring of
polynomials, and $I$ an almost complete intersection of finite
colength generated in degree $n$. Suppose that $I$ is $s$--balanced.
\begin{enumerate}
\item[{\rm (i)}] If $d=3$ and $s=n-1$, then $L_2$ is generated by ${{s+1}\choose{2}}$ forms with
coefficients of degree $s$ and
there are precisely $n$  linear syzygies of
degree $n-1$, and $n\leq 7\,${\rm ;}

\item[{\rm (ii)}] If $d=4$ and $s=n=2$, then there are precisely $15$ linear
syzygies of degree $2$.
\end{enumerate}
\end{Corollary}

\demo  We begin by observing the values of $r(I)$. In
 case (i), $r(I)= (3-1)(n-1)-s= (3-2)(n-1)=s$, while in (ii)
$r(I)= (4-1)(2-1)-2=1$.

\medskip

The first assertion of (i) comes from Proposition~\ref{syzlemma0} and
the value $r(I)=s$. As for the number of syzygies,
 according to Theorem~\ref{socdegree}(iv),
$R_{n+s}= I_{n+s}$ in case (i), and $R_{n+1}=I_{n+1}$ in case (ii),
which will permit the determination of the
dimension of the linear syzygies of degree $s$, or higher in case (i), and
for all degrees in case (ii).

\medskip

Let us focus on the case $r(I)=s$.
  Consider the exact sequence corresponding to the generators of
$I$,
\[ R^{d+1} \stackrel{\pi}{\lar} R \lar R/I\rar 0,\]
and write $\psi_s$ for the
vector space map
induced by $\pi$ on the homogeneous component
of degree $s$ of $R^{d+1}$. We have an exact sequence of $k$-vector spaces and $k$-linear
maps
\[ R_s^{d+1} \stackrel{\psi_s}{\lar} R_{n+s} \lar R_{n+s}/I_{n+s}\rar 0\]
and $\ker (\psi_s)$ is the $k$-span of the syzygies of $I$ of $R$-degree $s$.

One easily has
\begin{equation}\label{dimker} \dim_k (\ker (\psi_s))=
(d+1){{s+d-1}\choose{d-1}} - \dim_k(I_{n+s}).
\end{equation}

In this case, one gets
\begin{equation}\label{dimker2} \dim_k (\ker (\psi_s))=
(d+1){{s+d-1}\choose{d-1}} - {{s+n+d-1}\choose{d-1}}.
\end{equation}
If $I$ is moreover $s$-balanced for some $s\geq 1$ then it must be the case that
\begin{equation}\label{balanced_inequality} (d+1)\dim_k (\ker (\psi_s))\geq
\nu (\mathfrak{m}^s)=\dim_k(R_s)=
{{s+d-1}\choose{d-1}}.
\end{equation}

Suppose that $d=3$ and $I$ is $s$-balanced with $s=n-1$. Then the equality (\ref{dimker2}) gives
$ \dim_k (\ker (\psi_{n-1}))=n$ while the inequality (\ref{balanced_inequality})
easily yields $n\leq 7$.

\medskip

Finally, the assertion (ii) follows immediately from the equality (\ref{dimker2}).
\QED

\bigskip

The numerical data alone give a bird eye vision of the generators of the
graded pieces $L_1$ and $L_2$ of the ideal of equations of $L$.
This corollary is suitable in other cases, even when $I_1(\varphi)$ is a
less well packaged ideal.

\section{Binary Ideals}

In this section we take $d=2$ and write $R=k[x,y]$ (instead of the
general notation $R=k[x_1,x_2]$).
Let $I\subset R=k[x, y]$ be an $(x,y)$-primary ideal generated by three forms of
degree $n$. Suppose that  $I$ has a minimal free resolution \[ 0 \rar R^{2}
\stackrel{\varphi}{\lar} R^3 \lar I \rar 0.\]
We will assume throughout the section that the first   column of
$\varphi$ has  degree $r$, the other degree $s\geq r$. We note
$n=r+s$.

\medskip

Here we give a general format of the elimination equation of $I$ up to a power,
thus answering several questions raised in \cite{syl}.

\subsubsection*{Elimination equation and degree}

Set $S=R[\TT_1, \TT_2, \TT_3]$ as before.
Notice that $S$ is standard bigraded over $k$.
We denote by $f$ and $g$ the defining forms of the symmetric algebra
of $I$, i.e., the generators of the ideal $(L_1)\subset S$ in the earlier
notation. We write this in the form
\[ [f,g] = [\TT_1, \TT_2,\TT_3]\cdot  \varphi.\]
In the standard bigrading, by assumption, $f$ has bidegree
$(r,1)$, $g$ bidegree $(s,1)$.
According to Lemma~\ref{syzlemma}, the component $L_2$ could be determined from
$(L_1):I_1(\varphi)$. In dimension two it is more convenient to get
hold of a smaller quotient, $N=(L_1):(x,y)^{r}$. We apply basic
linkage theory to develop some properties of this ideal.

\begin{enumerate}
\item[$\bullet$]

 $N$, being a
direct link of the Cohen-Macaulay  ideal $(x,y)^r$, is a perfect
Cohen-Macaulay ideal of codimension two.
The canonical module of $S/N$ is generated by $(x,y)^rS/(f,g)$, so
that its Cohen-Macaulay type is $r+1$, according to
Theorem~\ref{Dubreil}.

\item[$\bullet$]  Therefore, by
the Hilbert-Burch theorem, $N$ is the ideal of maximal minors of an
$(r+2)\times (r+1)$ matrix $\zeta$ of homogeneous forms.

\item[$\bullet$] Thus,  $N=(f,g) : (x,y)^r$ has a
presentation $ 0 \rar S^{r+1} \stackrel{\zeta}{\lar} S^{r+2} \lar N
\rar 0$, where $\zeta$ can be written in the form \[ \zeta = \left[
\begin{array}{c} \sigma \\ ------- \\ \tau \end{array}\right], \]
with $\sigma$ is a $2 \times (r+1)$ submatrix with rows  whose entries are
biforms of bidegree  $(s-r, 1)$ and $(0,1)$; and $\tau$
is an $r \times (r+1)$ submatrix
whose entries are  biforms of bidegree $(1,0)$.
\item[$\bullet$] Since $N\subset L_2$, this shows that in $L_2$ there are
$r$ forms  $\hh_i$ of degree
$2$ in the $\TT_i$ whose $R$-coefficients are forms in $(x,y)^{s-1}$.

\item[$\bullet$]
If $s=r$, write
\begin{equation}\label{biforms}
[\hh]= [\hh_1 \;\; \cdots \;\;
\hh_r] = [x^{r-1}\;\; x^{r-2}y \;\; \cdots \;\; xy^{r-2} \;\; y^{r-1}]
\cdot \BB,
\end{equation}
 where $\BB$ is an $r \times r$ matrix whose entries belong
to $k[\TT_1, \TT_2, \TT_3].$

\item[$\bullet$]
If $s> r$, collect the $s-r$ forms $\ff=\{x^{s-r-1}f, x^{s-r-2}yf,\ldots, y^{s-r-1}f\}$
and write
\begin{equation}\label{biforms2}
[\ff;\hh]= [\ff; \hh_1 \;\; \cdots \;\;
\hh_r] = [x^{s-1}\;\; x^{s-2}y \;\; \cdots \;\; xy^{s-2} \;\; y^{s-1}]
\cdot \BB,
\end{equation}
 where $\BB$ is an $s \times s$ matrix whose entries belong
to $k[\TT_1, \TT_2, \TT_3]$.

\end{enumerate}

A first consequence of this analysis is one of our main results:

\begin{Theorem} \label{detB} In both cases,
 $\det \BB $ is a nonzero polynomial of degree $n$.
\end{Theorem}

\demo
Case $s=r$:
 Suppose that $\det \BB=0$. Then there exists a nonzero
vector $\left[\begin{array}{c} {\bf a}_1 \\ \vdots \\ {\bf a}_r
\end{array}\right]$  whose entries are in $k[T_1, T_2, T_3]$ such
that  $\BB \cdot \left[\begin{array}{c} {\bf a}_1 \\ \vdots \\ {\bf a}_r
\end{array}\right]=0$.  Hence $[\hh_1 \;\; \cdots \;\; \hh_r]
\left[\begin{array}{c} {\bf a}_1 \\ \vdots \\ {\bf a}_r  \end{array}\right]=0$.
Since the relations of  $ \hh_1,\ldots, \hh_r$ are $S$-linear combinations of
the columns of $\zeta$, we get a
contradiction.
\medskip

The assertion on the degree follows since the degree of $\det \BB$ is
$2r=r+s=n$.

\medskip

Case $s>r$:
 Suppose that $\det \BB=0$. Then there exists a nonzero
vector $\left[\begin{array}{c} {\bf a}_1 \\ \vdots \\ {\bf a}_s
\end{array}\right]$  whose entries are in $k[T_1, T_2, T_3]$ such
that  $\BB \cdot \left[\begin{array}{c} {\bf a}_1 \\ \vdots \\ {\bf
a}_s
\end{array}\right]=0$.  Hence $[\ff; \; \; \hh_1 \;\; \cdots \;\; \hh_r]
\left[\begin{array}{c} {\bf a}_1 \\ \vdots \\ {\bf a}_s
\end{array}\right]=0$.
We write this relation as follows
\[
\sum_{i=1}^{s-r} {\bf a_i}  x^{s-r-i} y^{i-1}f + \sum_{j=1}^r {\bf
a}_{s-r+j}\hh_j=
{\bf a}f
 + \sum_{j=1}^r {\bf
a}_{s-r+j}\hh_j=0,  \]
where \[ {\bf a}=
\sum_{i=1}^{s-r} {\bf a_i}  x^{s-r-i} y^{i-1}.\]

Since the relations of  $f, \hh_1,\ldots, \hh_r$ are $S$-linear combinations of
the columns of $\zeta$,
\begin{itemize}

\item ${\bf a}\in (x,y)^{s-r}S$ and
\item ${\bf a}_{s-r+j}\in (x,y)S$, for $ 1 \leq j \leq r $.
\end{itemize}
and therefore ${\bf a}_i\in (x,y)S$, for all $i$. This gives a
contradiction.
\medskip

The assertion on the degree follows since the degree of $\det \BB$ is
$(s-r) + 2r =r+s=n$.
\QED

\begin{Example}{\rm Let $R=k[x, y]$ and $I$ the ideal defined by
 {$\varphi=\left[ \begin{array}{ll}
x^2 & y^4 \\ xy & x^3y+x^4 \\ y^2 & xy^3
\end{array} \right]$.}
\begin{enumerate}

\item[$\bullet$] $N= (f,g) : \m^2 = (f,g, \hh_1, \hh_2)$, where
\[\begin{array}{lll}
f &= & x^2\TT_1+xy\TT_2+y^2\TT_3 \\
g &= & y^4\TT_1+(x^3y+x^4)\TT_2 + xy^3\TT_3 \\
\hh_1 &= &
y^3\TT_1^2-x^3\TT_2^2-x^2y\TT_2^2+xy^2\TT_1\TT_3-x^2y\TT_2\TT_3-xy^2\TT_2\TT_3 \\
\hh_2 &= &
xy^2\TT_1^2+y^3\TT_1\TT_2-x^3\TT_2\TT_3-x^2y\TT_2\TT_3-y^3\TT_3^2
\end{array}
\]

\item[$\bullet$] $[xf\;\; yf \;\; h_1 \;\; h_2 ] = \m^3 \BB$, where
\[
\BB=\left[
\begin{array}{llll}
\TT_1 & 0 & -\TT_2^2 & -\TT_2\TT_3 \\
\TT_2 & \TT_1 & -\TT_2^2-\TT_2\TT_3 & -\TT_2\TT_3 \\
\TT_3 & \TT_2 & \TT_1\TT_3 -\TT_2\TT_3 & \TT_1^2 \\
0 & \TT_3 & \TT_1^2 & \TT_1\TT_2-\TT_3^2
\end{array}
\right]
\]

\item[$\bullet$]  $\det\BB$ is the elimination equation.
\end{enumerate}
}\end{Example}

\begin{Remark}{\rm The polynomials $\hh_1, \ldots, \hh_r$ were also obtained
in \cite{syl} by a direct process involving Sylvester elimination, in
the cases $s=r$ or $s=r+1$. In
\cite{syl} though they did not arrive with the elements of
structure--that is with their relations--provided in the
Hilbert-Burch matrix. It is this fact that opens the way in the
binary case to a greater generality to the ideals treated and a more
detailed understanding of the ideal $L$.
}\end{Remark}

The next result provides a secondary elimination degree for these
ideals.

\begin{Corollary}\label{balanced_exponent} $L= (L_1): (x,y)^{n-1}$.
\end{Corollary}
\demo
With the previous notation, let $N=L\cap Q$ be the primary decomposition of $N$,
where $Q$ is $(x,y)S$-primary. Writing $\beta=\det \BB$, we then have $N:\beta=Q$.
On the other hand, $(L_1,{\bf h})\subset N$ by construction and
$(x,y)^{s-1}S\subset (L_1,{\bf h}):\beta$
since by (\ref{biforms} and \ref{biforms2}) the biforms $\hh$, or
$\ff, \hh$,  must effectively involve all monomials of degree $s-1$ in $x,y$.
It follows that $(x,y)^{s-1}S\subset Q$, hence
$$L(x,y)^{s-1}\subset LQ\subset L\cap Q=N=(L_1):(x,y)^r,$$
thus implying that $L(x,y)^{s+r-1}\subset (L_1)$.
This shows the assertion.
\QED

\subsubsection*{Elimination equation up to a power}

\begin{Theorem} \label{detB2plus} Let $I$ be as above and $\beta=\det \BB$. Then
$\beta$ is a power of the elimination equation of $I$.
\end{Theorem}

\demo
Let $\pp$ denote the elimination equation of $I$.
Since $\pp$ is irreducible it suffices to show that $\beta$ divides
a power of $\pp$.

The associated primes of $N=(L_1): (x,y)^r$ are
the defining ideal $L$ of the Rees algebra and
 $\mathfrak{m}S= (x,y)S$.
 We have a primary decomposition
\[ N= L \cap Q,\]
where $Q$ is $\mathfrak{m}S$-primary.
  From the proof of
  Theorem~\ref{detB}, localizing at $\mathfrak{m}S$ gives
$ (x,y)^{s-1}S\subset Q$. (Equality will hold
  when $r=s$.)

 The equality $N= L\cap Q= (L_1, \hh) $ implies that
$(x,y)^{s-1}\pp \subset (L_1,\hh)$.
Since $f,g$ are of bidegrees $(r,1)$ and $(s,1)$, it must be
the case that each polynomial
 $x^{i}y^{s-1-i}\pp $ lies in the span of the $(f(x,y)^{s-r-1}, \hh)$
alone.   This gives a representation
\[ \pp[(x,y)^{s-1}]= [\ff;\hh] \cdot \AA,\]
(or simply $\pp[(x,y)^{r-1}]= [\hh] \cdot \AA$, if $s=r$)
where $\AA$ is an $s\times s$ matrix with entries in $k[\TT_1, \TT_2, \TT_3]$.
Replacing $[\ff;\hh]$ by $[(x,y)^{s-1}] \cdot \BB$, gives the matrix
equation
\[ [(x,y)^{s-1}] \big( \BB\cdot \AA-\pp \mathbf{I}\big)=0, \]
where $\mathbf{I}$ is the $s\times s$ identity matrix.

Since the minimal syzygies of $(x,y)^{s-1}$ have coefficients in
$(x,y)$, we must have
 \[\BB\cdot \AA=\pp \mathbf{I}, \]
so that $\det \BB\cdot \det \AA= \pp^r$, as desired. \QED


\section{Ternary Ideals}\label{ternary}

We outline a conjectural scenario that we expect many such
ideals  to conform to. Suppose $I$ is an ideal of $R=k[x_1,x_2,x_3]$
generated by forms
$a_1,a_2, a_3, a$ of degree $n\geq 2$, with $J=(a_1,a_2,a_3)$ being a
minimal reduction. This approach is required because the linkage
theory method lacks the predicability  of the  binary ideal case.

\bigskip

\subsubsection*{Balanced ternary ideals}

Suppose that $I$ is $(n-1)$--balanced, where $n$ is the degree
 of the generators of $I$.
By Corollary~\ref{wellbal}(1),
 there are $n$ linear forms $\ff_i$, $1\leq
i\leq n$,
\[ \ff_i= \sum_{j=1}^4 c_{ij}\TT_j\in L_1,\]
arising  from the syzygies of $I$ of degree $n-1$. These syzygies,
according to Theorem~\ref{Dubreil}, come from the syzygies of $J:a$,
which by the structure theorem of codimension three Gorenstein
ideals, is given by the Pfaffians of a skew-symmetric matrix $\Phi$,
of size  at most $2n-1$.

According to
Proposition~\ref{syzlemma0}, there are
$n\choose 2$ quadratic forms $\hh_k$  ($1\leq k\leq {{n}\choose
{2}}$):
\[ \hh_k= \sum_{1\leq i\leq j\leq 4} c_{ijk}\TT_i\TT_j\in L_2,\]
with $R$-coefficients of degree $n-1$.

Picking a basis for $\mathfrak{m}^{n-1}$ (simply denoted by
$\mathfrak{m}^{n-1}$),
 and writing the $\ff_i$ and $\hh_k$ in matrix format, we have
\begin{eqnarray} \label{acs3sedeg}
 [\ff_1, \ldots, \ff_n, \hh_1, \ldots, \hh_{{n}\choose{2}}]&=&
\mathfrak{m}^{n-1}\cdot
\BB,
\end{eqnarray}
where $\BB$ is the corresponding content matrix (see \cite{syl}).
Observe that $\det\BB$ is either zero, or a polynomial of degree
\[ n+ 2{{n}\choose{2}}= n^2.\]
It is therefore a likely candidate for the {\em elimination
equation}. Verification
consists in checking that $\det\BB$ is irreducible for an ideal in
any
given {\em generic} class. We will make this more precise
 on a quick analysis of the lower degree cases.

\begin{Theorem} \label{detB3s} If $I\subset R=k[x_1,x_2,x_3]$ is a  $(n-1)$-balanced almost complete
intersection ideal  generated by forms of degree $n$ {\rm (}$n\leq 7${\rm )}, then
\begin{eqnarray*}
 \det \BB\neq 0.
\end{eqnarray*}
\end{Theorem}

\demo  Write each of the
quadrics ${\hh}_j$ above in the form
\[ {\hh}_j = c_j\TT_4^2 + \TT_4\ff(\TT_1,\TT_2, \TT_3)+
\g2(\TT_1,\TT_2, \TT_3),
 \]
where $c_j$ is a form of $R$ of degree $n-1$, and similarly,
\[ \ff_i= c_{i4} \TT_4+ \sum_{k=1}^3 c_{ik}\TT_k \in (L_1).\]

\medskip

Write $\mathfrak{c}=(c_j, c_{i4})$ for
the ideal of $R$ generated by the leading coefficients of
$\TT_4$ in the $\ff_i$'s and of $\TT_4^2$ in the ${\hh}_j$'s. It is apparent that if $\mathfrak{c}=\mathfrak{m}^{n-1}$, there will be a non-cancelling
term $\TT_4^{n^2}$ in the expansion of $\det \BB$.

To argue that indeed $\mathfrak{c}=\mathfrak{m}^{n-1}$ is the case, assume otherwise.
Since the $\ff_i$ are minimal generators
that contribute to $(J:a)$, we may assume that the $c_{i4}$ are
linearly independent. This implies that  we may replace one of the ${\hh}_j$ by a
form
\begin{eqnarray*}{\hh} &=&  \TT_4\ff(\TT_1,\TT_2, \TT_3)+
\g2(\TT_1, \TT_2, \TT_3)\\
&=&
  \TT_4(r_1\TT_1+  r_2\TT_2+r_3 \TT_3)+
\TT_1\g2_1+  \TT_2\g2_2+ \TT_3\g2_3,
\end{eqnarray*}
where the $r_i$ are $(n-1)$-forms in $R$ and the $\g2$'s are $\TT$-linear
involving  only $\TT_1,\TT_2,\TT_3$.

Evaluate now $\TT_i$ at the corresponding generator of $I$ to get \[
(ar_1+\g2_1(a_1,a_2,a_3) )\, a_1 + (ar_2+\g2_2(a_1,a_2,a_3) )\, a_3+
(ar_3+\g2_3(a_1,a_2,a_3) )\, a_3=0,\] a syzygy of the ideal $J=(a_1,a_2,a_3)$.
Since $J$ is a complete intersection, $ar_i+\g2_i(a_1,a_2,a_3)={\bf
u}_i(a_1,a_2,a_3)\in J$ for $i=1,2,3$, with ${\bf u}_i$ a linear form
in $\TT_1,\TT_2,\TT_3$ with coefficients in $R$.  These are syzygies
of the generators of $I$, so lifting back to $1$-forms in $\TT$ and
substituting yields $\hh=\hh'+{\bf k}$, where
\begin{eqnarray*} \hh'&=&(r_1\TT_4+\g2_1(\TT_2,\TT_2,\TT_3)-{\bf u}_1(\TT_1,\TT_2,\TT_3))\,\TT_1\\
&+& (r_2\TT_4+\g2_2(\TT_1,\TT_2,\TT_3)-{\bf u}_2(\TT_1,\TT_2,\TT_3))\,\TT_2\\
&+& (r_3\TT_4+\g2_3(\TT_1,\TT_2,\TT_3)-{\bf u}_3(\TT_1,\TT_2,\TT_3))\,\TT_3
\end{eqnarray*}
is an element of $(L_1)$, and because $\hh$ is a relation then so is the term ${\bf k}$.
But the latter only involves $\TT_1,\TT_2,\TT_3$, hence it belongs to the defining ideal of the symmetric
algebra of the complete intersection $J$. But the latter is certainly contained in $(L_1)$.
Summing up we have found that $\hh\in (L_1)$, which is a contradiction
since $\hh$ is a minimal
generator of $L_2/S_1L_1$. \QED

\begin{Example} {\rm This example shows, in particular, that there
are $(n-1)$--balanced ideals $I\subset k[x_1,x_2,x_3]$ generated in degree $n$
such that the corresponding map $\Psi_I$ is not birational onto its image.

 Let  $I=(J,a)$, where
$J=(x_1^3,\,x_2^3,\,x_3^3)$ and $a=x_1x_2x_3$.

\begin{enumerate}

\item[$\bullet$] The Hilbert series of $R/(J:a)$ is $1+3t+3t^2+t^3$.

\item[$\bullet$] $I_1(\varphi)=\m^2$, i.e., $I$ is $2$--balanced.

\item[$\bullet$] $(L_1) : \m^4 = (L_1) : \m^5$ while $(L_1) : \m^3 \neq (L_1) : \m^4$

\item[$\bullet$] Equations of $I$:
\[\left\{\begin{array}{ll} L=(L_1,L_2,L_3)& \\
 \nu(L_1)=6; \nu(L_2)=3; L_3=kF, &
\end{array}
\right.
\]
where $F=-\TT_1\TT_2\TT_3+\TT_4^3$ is the elimination equation; in particular, the
corresponding map $\Psi_I$ is not birational onto its image.

\item[$\bullet$] Let $\ff_1, \ff_2, \ff_3$ be generators of
$L_1$ with coefficients in $\m^2$ and $\hh_1, \hh_2, \hh_3\in L_2$
as previously described. Writing $[\ff_1, \ff_2, \ff_3, \hh_1, \hh_2, \hh_3]=\m^2 \BB$
as in (\ref{acs3sedeg}), one has
\[ \BB = \left[\begin{array}{rrrrrr} 0& 0&
-\TT_4 & \TT_3 & 0 & 0 \\ 0 & -\TT_4 & 0 & 0 & \TT_2 & 0 \\ -\TT_4 & 0& 0 & 0 &
0& \TT_1 \\ \TT_2 \TT_3 & 0 & 0& 0& 0 &-\TT_4^2 \\ 0 & \TT_1T_3 &0&0
&-\TT_4^2 &0
\\ 0 & 0 & \TT_1\TT_2 &-\TT_4^2 &0 & 0 \end{array} \right]
\]
and $\det \BB = F^3$.
\end{enumerate}
}\end{Example}

\begin{Remark}{\rm To strengthen Theorem~\ref{detB3s} to the
assertion that $\det \BB$ is a power of the elimination equation, one
needs more understanding of the ideal $(L_1):\mathfrak{m}^{n-1}$.
Here is one such instance.
}\end{Remark}

\begin{Proposition}\label{being_a_power}
Let $I\subset R=k[x_1,x_2,x_3]$ be an $(n-1)$-balanced almost complete
intersection ideal  generated by forms of degree $n$.
Keeping the notation introduced at the
beginning of this section, we obtain the following :
\begin{enumerate}
\item[{\rm (i)}] $2n-2$ is a secondary elimination degree of $I$.

\item[{\rm (ii)}] If $(L_1):\mathfrak{m}^{n-1} = ({\bf f}, {\bf h})=(\ff_1, \ldots, \ff_n, \hh_1, \ldots, \hh_{{n}\choose{2}})$, then  $\det \BB$ is a power of the elimination equation.
\end{enumerate}
\end{Proposition}

\demo
Let $(L_1):\mathfrak{m}^{n-1}=L\cap Q$, where $Q$ an
$\mathfrak{m}S$-primary ideal. As in  (\ref{acs3sedeg}) one has
$ [\ff, \hh]= \m^{n-1}\cdot \BB$. Notice that  $(\ff, \hh) \subset ( (L_1) : \m^{n-1}) $.
Write $\beta=\det \BB$. Since $\det \BB \neq 0$,
it follows that
\[ \m^{n-1}S \subset (\ff, \hh) : \beta \subset ( (L_1) : \m^{n-1}) : \beta =Q.\]
This implies that
\[ L \cdot \mathfrak{m}^{n-1} \subset  LQ
\subset (L_1):\mathfrak{m}^{n-1}.\]
Hence
$L= (L_1):\mathfrak{m}^{2n-2}$, which proves (i).

Now suppose that $(L_1):\mathfrak{m}^{n-1} = ({\bf f}, {\bf h})$. Let $\pp$ be the elimination equation.  By (i), we have
 \[\pp \in (L_1):\mathfrak{m}^{2n-2} = ( (L_1):\m^{n-1} ) : \m^{n-1}. \] Therefore
$\pp \mathfrak{m}^{n-1}\subset (\ff, \hh)$, which gives a representation
\[ \pp [\mathfrak{m}^{n-1}]= [\ff, \hh] \cdot \AA,\]
where $\AA$ is a square matrix with entries in $S$.
Replacing $[\ff,\hh]$ by $[\mathfrak{m}^{n-1}] \cdot \BB$, gives the matrix
equation
\[ [\mathfrak{m}^{n-1}] \big( \BB\cdot \AA-\pp \mathbf{I}\big)=0, \]
where $\mathbf{I}$ is the  identity matrix. Since the minimal syzygies of $\mathfrak{m}^{r-1}$ have coefficients in
$\mathfrak{m}$, we must have
 \[\BB\cdot \AA=\pp \mathbf{I}, \]
so that $\det \BB\cdot \det \AA= \pp^m$, for some integer $m$; this proves (ii).
 \QED

\subsubsection*{Ternary quadrics}

We apply the preceding discussion to the situation where the ideal
$I$ is generated by $4$ quadrics of the
polynomial ring $R=k[x_1,x_2,x_3]$. The socle of
$R/J$ is generated by the Jacobian determinant of $a_1,a_2,a_3$,
which implies that $\lambda(I/J)\geq 2$. Together we obtain that
$\lambda(R/I)=6$. The Hilbert function of $R/I$ is $(1,3, 2)$. Since
we cannot have $u\mathfrak{m}\subset I$ for some $1$-form $u$, the
type of $I$ is $2$ and its socle is generated in degree two.

The canonical module of $R/I$ satisfies  $\lambda((J:a)/J)=6$, hence
$\lambda(\mathfrak{m}/(J:a))=1$, that is to say
$J:a=(v_1,v_2,v_3^2)$, where the $v_i$ are linearly independent
$1$-forms. Let $\ff_1$ and $\ff_2$ be the linear syzygies of $I$ induced by $v_1$ and $v_2$ respectively.
$R/I$ has a free presentation
\[ 0 \rar R^2\lar R^5\stackrel{\varphi}{\lar} R^4\lar R/I \rar 0.\]

The ideal $I_1(\varphi)$ is either
$\mathfrak{m}$ or  $(v_1,v_2,v_3^2)$.
 In the first case, $\delta(I)$ is the socle of
$H_1(I)$, an element of degree $3$,
 therefore its image in $ {L}_2$ is a $2$-form $\hh_1$
with linear coefficients. In particular reduction number cannot be
two.
Putting it together with the two
linear syzygies $\ff_1, \ff_2$ of $I$, we
\[ [\ff_1, \ff_2, \hh_1] = [x_1,x_2,x_3]\cdot \BB, \]
where $\BB$ is a $3\times 3$ matrix with entries in $k[\TT_1,
\TT_2,\TT_3,\TT_4]$, of column degrees $(1,1,2)$. The quartic $\det \BB$,
 is the elimination equation of $I$.

In the other case, $\delta(I)$ has degree two, so its image in
$ {L}_2$ is a form with coefficients in the field. The
corresponding mapping  $\Psi_I$ is not birational.

\medskip

We sum up the findings in this case:

\begin{Theorem} \label{aci32} Let $R=k[x_1, x_2, x_3]$ and let $I$ be an
$\mathfrak{m}$-primary
almost complete intersection generated by quadrics. Then
\begin{enumerate}
\item[{\rm (i)}]  If $I_1(\varphi)=\mathfrak{m}$ the corresponding
mapping $\Psi_I$ is
birational onto its image.
\item[{\rm (ii)}]  If  $I_1(\varphi)\neq\mathfrak{m}$ then $I_1(\varphi)= (v_1, v_2, v_3^2)$, where $v_1,v_2,
v_3$ are linearly independent $1$-forms,  and the mapping $\Psi_I$ is
not birational onto its image.
\end{enumerate}
\end{Theorem}

Here is a sufficiently general example fitting the first case in the above theorem
(same behavior as $4$ random quadrics):
\begin{equation}\label{quadrics}
J=(x_1^2, x_2^2, x_3^2), \quad a=x_1x_2+x_1x_3+x_2x_3.
\end{equation}
A calculation shows that $\det \BB$  is irreducible.

\bigskip

An example which is degenerate (non-birational) as in the second case above is

\begin{equation}\label{quadrics_degenerate}
J=(x_1^2, x_2^2, x_3^2), \quad a=x_1x_2.
\end{equation}
Then $(L_1, x_3(\TT_1 \TT_2-\TT_4^2)) = (L_1): \m \subsetneq (L_1): I_1(\varphi)=L$.
Let $\hh=x_3(\TT_1 \TT_2-\TT_4^2)$ and write \[ [\ff_1, \ff_2, \hh] = [x_1,x_2,x_3]\cdot \BB. \]
Then $\det \BB $ is a square of the elimination equation $\pp=\TT_1 \TT_2-\TT_4^2$, so we still recover the
elimination equation from $\BB$.

\subsubsection*{Ternary cubics and quartics}

Let $I$ be an ideal generated by $4$ cubics and suppose that $I$ is
$2$-balanced (i.e.,
$I_1(\varphi)=\mathfrak{m}^2$).
Using this (see the beginning of Section~\ref{ternary}) and the fact that
 $J:a$ is a codimension $3$ Gorenstein ideal, it follows that $J:a$ is minimally generated by the Pfaffians of
a skew-symmetric matrix of sizes $3$ or $5$.

In the first case, $J:a$
is generated by $3$ quadrics and $I$ is a Northcott ideal. In the
second  case, $J:a$
cannot be generated by $5$ quadrics, as its Hilbert function would be
$(1,3,1)$ and therefore the Hilbert function of $R/I$ would have to be
\[ (1,3,6,7,6,3,1)- (0,0,0, 1,3,1)= (1,3,6,6, 3,2,1),\]
giving that the canonical module of $R/I$ had a generator in degree
$0$. Thus $J:a$ must be generated by
 $3$ quadrics and $2$ cubics.

\begin{Remark} \label{aci33and4}{\rm  Let $R=k[x_1, x_2, x_3]$ and $I$ a
  balanced
$\mathfrak{m}$-primary
almost complete intersection. We expect that with an appropriate
notion of genericity the following assertions will hold.
\begin{enumerate}
\item If $I$
 generated by cubics and  $I_1(\varphi)=
\mathfrak{m}^2$,
 the polynomial
$\det \BB$,
defined by the equation {\rm (\ref{acs3sedeg})}, is the elimination
equation of $I$.

\item If $I$
 generated by quartics and  $I_1(\varphi)=
\mathfrak{m}^3$,
 the polynomial
$\det \BB$,
defined by the equation {\rm (\ref{acs3sedeg})}, is the elimination
equation of $I$.
\end{enumerate}
}\end{Remark}

In what follows we give examples to cover this expected behavior:
the first two are instances of (1), while the third illustrates (2).

\begin{Example}
{\rm Let  $I=(J, a)$, where
$J=(x_1^3+x_2^2x_3,\,x_2^3+x_1x_3^2,\, x_3^3+x_1^2x_2)$ and $a=x_1x_2x_3$.

\begin{enumerate}

\item[$\bullet$] $J:a$ is a complete intersection



\item[$\bullet$] $(L_1) : \m^4 = (L_1) : \m^5$ while $(L_1) : \m^3 \neq (L_1) : \m^4$

\item[$\bullet$] Equations of $I$:
\[\left\{\begin{array}{ll} L=(L_1,L_2,L_5,L_9)& \\
 \nu(L_1)=6; \,\nu(L_2)=3; \,\nu(L_5)=15;\, L_9=k\det \BB, &
\end{array}
\right.
\]
where $\det \BB$, obtained as in (\ref{acs3sedeg}), is of degree $9$, hence must be the
elimination equation and
the mapping $\Psi_I$ is birational onto its image.


\end{enumerate}
}
\end{Example}

\begin{Example}
{\rm Let $I=(J, a)$, where
$J=(x_1^3,\,x_2^3,\,x_3^3)$ and $a=x_1^2x_2+x_2^2x_3+x_1x_3^2$.

\begin{enumerate}
\item[$\bullet$] $J:a$ is Gorenstein.



\item[$\bullet$] $(L_1) : \m^4 = (L_1) : \m^5$ while $(L_1) : \m^3 \neq (L_1) : \m^4$.

\item[$\bullet$] Equations of $I$:
$$\left\{\begin{array}{ll} L=(L_1,L_2,L_4,L_9)& \\
 \nu(L_1)=7;\, \nu(L_2)=3;\, \nu(L_4)=6;\, L_9=k\det \BB, &
\end{array}
\right.
$$
where $\det \BB$, obtained as in (\ref{acs3sedeg}), is of degree $9$, hence must be the elimination equation and
the mapping $\Psi_I$ is birational onto its image.


\end{enumerate}

}
\end{Example}

\begin{Example}{\rm
Let $I=(J, a)$, where
$J=(x_1^4,\,x_2^4,\,x_3^4)$ and $a=x_1^3x_2+x_2^3x_3+x_1x_3^3$.

This follows the pattern of the previous example, with
 $J:a$ is Gorenstein.
 The structure of $L$ is now involved, however the principle in
 (\ref{acs3sedeg}) still works and gives
$\det \BB$ of degree $16$, hence must be the elimination equation and
the mapping $\Psi_I$ is birational onto its image.
}
\end{Example}

\section{Quaternary Forms}

In this Section, we set $R=k[x_1, x_2,x_3, x_4]$ with $\m=(x_1, x_2,x_3, x_4)$.

\subsubsection*{Quaternary quadrics}

Let $I$ be generated by $5$ quadrics $a_1, a_2, a_3, a_4, a_5 $ of
$R$,  $J=(a_1, a_2, a_3, a_4)$, and $a=a_5$.
The analysis of this case is
less extensive than the case of ternary ideals. Let us assume that
$I_1(\varphi)=\mathfrak{m}^2$ -- that is,  the balancedness exponent equals the degree
of the generators, a case that occurs generically in this degree.

We claim that the Hilbert function of $R/I$ is $(1,4,5)$.
Since
we cannot have $u\mathfrak{m}\subset I$ for some $1$-form $u$,
 its socle is generated in degree two or higher. First, we argue that
 $J:a\neq \mathfrak{m}^2$; in fact,  otherwise $R/I$  would be of length $11$
 and type $6$ as $\mathfrak{m}^2/J$
 is its canonical module. But then the Hilbert function of  $R/I$ would be $(1,4,5,1)$, and the
 last two graded components would be in the  socle, which is impossible.

Thus it must be the case that $\lambda(R/I)=10$ and the Hilbert function of $H_1(I)$ is
$(5,4,1)$. Note that $\nu(\delta(I))=4$.

The last two graded components are of degrees $3$ and $4$. In degree
$3$ it leads to $4$ forms $\qq_1, \qq_2, \qq_3, \qq_4$ in $L_2$, with linear
coefficients:
\[ [\qq_1, \qq_2, \qq_3, \qq_4]= [x_1, x_2, x_3, x_4]\cdot \BB.\]

\begin{Theorem} \label{detB42}
$\det \BB\neq 0$.
\end{Theorem}

\demo We follow the pattern of the proof of Theorem~\ref{detB3s}.  Write each of the
quadrics $\qq_i$ as
\[ \qq_i = c_i\TT_5^2 + \TT_5\ff(\TT_1 \ldots, \TT_4)+
\g2(\TT_1, \ldots, \TT_4),
 \]
where $c_i$ is a linear form of $R$. If $(c_1,c_2,c_3, c_4)= \mathfrak{m}$
we  can take the $c_i$ for indeterminates in order to obtain the
corresponding $\det \BB$. In this case the occurrence of a non-cancelling term
$\TT_5^8$ in $\det\BB$ would be clear.

By contradiction, assume that the forms $c_i$ are not linearly
independent. In this case, we could replace one of the $\qq_i$ by a
form \[ \qq=  \TT_5\ff(\TT_1, \ldots, \TT_4)+
\g2(\TT_1, \ldots, \TT_4). \]
Keeping in mind that $\qq$ is a minimal generator we are going to argue
that $\qq \in (L_1)$. For that end, write
the form  $\qq$ as
\[ \qq=
  \TT_5(r_1\TT_1+  \cdots+r_4 \TT_4)+
\TT_1\g2_1+  \cdots + \TT_4\g2_4,  \]
where the $r_i$ are $1$-forms in $R$ and $\g2_i$'s are $\TT$--linear involving only $\TT_1, \ldots, \TT_4$.

Evaluate now the leading $\TT_i$ at the ideal to get the syzygy
\[ \hh= a(r_1\TT_1+  \cdots+r_4 \TT_4)+
a_1\g2_1+  \cdots + a_4\g2_4,  \]
of $I$, but actually of the ideal $J$. Since $J$ is a complete
intersection, all the coefficients of $\hh$ lie in $J$.
This implies that $r_ia\in J$ for all $r_i$. But this is impossible
since $J:a\in I_1(\varphi)=\mathfrak{m}^2$, unless all $r_i=0$.
This would imply that $\qq\in (L_1)$, as asserted. \QED

\begin{Theorem} \label{aci42} Let $R=k[x_1, x_2, x_3, x_4]$ and $I$ an
$\mathfrak{m}$-primary almost complete intersection.
 If $I$ is generated by quadrics and  $I_1(\varphi)=
\mathfrak{m}^2$, the polynomial $\det \BB$,
defined by the equation {\rm (\ref{acs3sedeg})} is divisible by the
elimination equation of $I$.
\end{Theorem}

\subsubsection*{Primary decomposition}

We now derive a value for the secondary elimination degree via a
primary decomposition.

\begin{Proposition} Let $R=k[x_1, x_2, x_3, x_4]$ and $I$ an
$\mathfrak{m}$-primary almost complete intersection.
Suppose that $I$ is generated by quadrics and  $I_1(\varphi)=
\mathfrak{m}^2$.
Then \[ (L_1): \mathfrak{m}^2= L \cap \mathfrak{m}S.\]
In particular,
\[L= (L_1): \mathfrak{m}^3. \]
\end{Proposition}

By Theorem~\ref{secondedeg}, it will suffice to show:

\begin{Lemma} $\mathfrak{m}^3= (J:a)\mathfrak{m}$.
\end{Lemma}

\demo We make use of the diagram

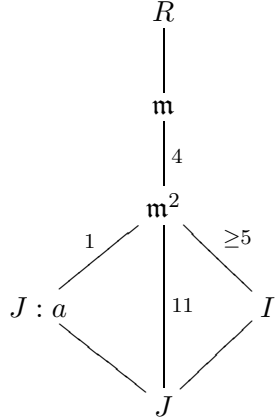
\[\diagram
& R \dline
  & \\
& \mathfrak{m}\dline^{4} & \\
& \mathfrak{m}^2\dlline_{1} \drline^{\geq 5} \ddline^{11} & \\
 J:a \drline & & I \dlline \\
 & J & \\
\enddiagram
\]

Let $I=(a_1, \ldots, a_5)$,  $J=(a_1, \ldots, a_4)$ a
minimal reduction of $I$, and $a=a_5$.
Since $(J:a)$ is a Gorenstein ideal, it cannot be $\m^2$.
On the other hand,
$\lambda((J:a)/J)= \lambda(R/I)\geq 10$. Thus
$\lambda(\mathfrak{m}^2/(J:a))=1$ and
$(J:a):\mathfrak{m}=\mathfrak{m}^2$ defines the socle of $R/(J:a)$.
The Hilbert function of $R/(J:a)$ is then $(1,4,1)$ which implies that
$\mathfrak{m}^3= (J:a)\mathfrak{m}$. \QED

\bigskip

\subsubsection*{Examples of quaternary quadrics}

We give a glimpse of the various cases.

\begin{Example}{\rm Let $J=(x_1^2,\,x_2^2,\,x_3^2,\, x_4^2)$
and $a=x_1x_2+x_2x_3+x_3x_4$

\begin{enumerate}
\item[$\bullet$] $\nu(J:a)=9$ (all quadrics)

\item[$\bullet$] $I_1(\varphi)=\mathfrak{m}^2$


\item[$\bullet$] As explained above, there are
four Rees equations $q_j$ of bidegree $(1,2)$ coming from the syzygetic principle.
One has:
\[\left\{\begin{array}{ll}
L=(L_1,L_2,L_9),& \\
\nu(L_1)=15; \nu(L_2)=4;  L_9=k\pp, & 
\end{array}
\right.
\]
where $\pp=\det\BB$ has degree $8$, hence is the elimination equation and the
corresponding map is birational.
\end{enumerate}
}
\end{Example}

The following example is similar to the above example, including the syzygetic principle, except that
$F$ is now the square root of $\det\BB$.

\begin{Example}{\rm Let
$J=(x_1^2,\,x_2^2,\,x_3^2,\, x_4^2)$ and  $a=x_1x_2+x_3x_4.$

\begin{enumerate}
\item[$\bullet$] $\nu(J:a)=9$ (all quadrics)

\item[$\bullet$] $I_1(\varphi)=\m^2$



\item[$\bullet$] Equations of $I$:
\[\left\{\begin{array}{ll} L=(L_1,L_2,L_4)& \\
 \nu(L_1)=15; \nu(L_2)=4; L_4=k\pp, &
\end{array}
\right.
\]
where $\deg(\pp)=4$, hence $\Psi_I$ is not birational onto its image.

\item[$\bullet$] Letting $\qq_1, \qq_2, \qq_3, \qq_4$  be generating forms in $L_2$, with linear
coefficients, write
\[ [\qq_1, \qq_2, \qq_3, \qq_4]= [x_1, x_2, x_3, x_4]\cdot \BB.\]
Then $\det \BB =\pp^2$.
\end{enumerate}
}\end{Example}

Next is an example where the normal syzygetic procedure fails, but one can
apply one more step to get the elimination equation.
We will accordingly give the details of the calculation.

\begin{Example}{\rm Let $I=(J,a)$, where
$J=(x_1^2,\,x_2^2,\,x_3^2,\, x_4^2)$ and $a=x_1x_2+x_2x_3+x_3x_4+x_1x_4.$

\begin{enumerate}
\item[$\bullet$] $\nu(J:a)=4$ (complete intersection of $2$ linear equations and $2$ quadrics)

\item[$\bullet$] $I_1(\varphi)=\m$ (i.e., $1$-balanced, not $2$-balanced)

\item[$\bullet$] $\nu(Z_1)=9$ ($2$ linear syzygies and $7$ quadratic ones)

\item[$\bullet$] Hilbert series of $R/(J:a)$ : $1+2t+t^2$.


\item[$\bullet$] $\edeg(I)=8=2^3$ (i.e., birational)

\item[$\bullet$] The one step syzygetic principle does not work to get the
elimination equation:
One has:
\[\left\{\begin{array}{ll}
L=(L_1, L_2, L_3, L_8),& \\
\nu(L_1)=9; \nu(L_2)=1;  \nu(L_3) =2, L_8=k\pp, &
\end{array}
\right.\]
Nevertheless there are relationships among these numbers that are
understood from the syzygetic discussion. Thus
Theorem~\ref{socdegree}(iii) implies that $L_2/S_1L_1$ is generated
by one form with linear coefficients.

\item[$\bullet$] Let $f_1$ and $f_2$ be in $L_1$ with linear coefficients, i.e.,
\begin{enumerate}
\item[$\diamondsuit$] $f_1=(x_1+x_3)(\TT_2-\TT_4)+(x_2-x_4)(-\TT_5)$

\item[$\diamondsuit$] $f_2=(x_1-x_3)(-\TT_5)+(x_2+x_4)(\TT_1-\TT_3)$
\end{enumerate}

\item[$\bullet$] Let $\hh+L_1S_1$ be a generator of $L_2/L_1S_1$. We
observed that, as a coset,  $\hh$ can be written as $\hh=\hh_1+S_1L_1$ and
$\hh=\hh_2+S_1L_1$, with $\hh_1$ and $\hh_2$ forms of bidegree $(2,2)$, with
$R$-content contained in the contents $f_1$ and $f_2$, respectively.

\item[$\bullet$] Let
$\hh_1=(-2x_1\TT_4\TT_5-2x_3\TT_4\TT_5+x_4\TT_5^2)(x_1+x_3)+(x_3\TT_1\TT_2-x_3\TT_2\TT_3-x_3\TT_1\TT_4+x_3\TT_3\TT_4-
x_4\TT_1\TT_5+2x_2\TT_3\TT_5-x_4\TT_3\TT_5-x_3\TT_5^2)(x_2-x_4)$. Then $\hh_1 \in
L_2$ and $2\hh+\hh_1 \in L_1S_1$.  Hence we may choose $\hh_1+L_1S_1$ to be
a generator of $L_2/L_1S_1$.

\item[$\bullet$] Let $\hh_2=(x_4\TT_1\TT_2-x_4\TT_2\TT_3-x_4\TT_1\TT_4+x_4\TT_3\TT_4-x_3\TT_2\TT_5+2x_1\TT_4\TT_5-x_3\TT_4\TT_5-x_4\TT_5^2)(x_1-x_3)+
(-2x_2\TT_3\TT_5-2x_4\TT_3\TT_5+x_3\TT_5^2)(x_2+x_4)$. Then $\hh_2 \in L_2$ and
$2\hh+\hh_2 \in L_1S_1$.
Hence we may choose $\hh_2+L_1S_1$ to be a generator of $L_2/L_1S_1$.

\item[$\bullet$] Write
\[
[ f_1 \;\; \hh_1]  =  [x_1+x_3 \;\; x_2-x_4]\BB_1 \quad \mbox{\rm and} \quad
[ f_2 \;\; \hh_2]  = [x_1-x_3 \;\; x_2+x_4] \BB_2
\]
Then $\det \BB_1$ and $\det \BB_2$ form a minimal generating set of
$L_3/S_1L_2$.

\item[$\bullet$] Write $[f_1\;\; f_2\;\; \det \BB_1 \;\; \det \BB_2]
= \  [x_1 \;\; x_2 \;\; x_3 \;\; x_4]\BB$.

At the outcome
$\det \BB = \pp$ is of degree $8$, hence this is again birational.
\end{enumerate}

}\end{Example}

\end{document}